\documentclass[10pt]{amsart}
     \makeatletter
     \def\section{\@startsection{section}{1}%
     \z@{.7\linespacing\@plus\linespacing}{.5\linespacing}%
     {\bfseries
     \centering
     }}
     \def\@secnumfont{\bfseries}
     \makeatother
\newtheorem{theorem}{Theorem}[section]

\theoremstyle{definition}

\theoremstyle{remark}

\numberwithin{equation}{section}
\usepackage{amsmath,amsthm,amssymb,amsbsy}

\newcommand{\QQ}{\mathbb{Q}}
\newcommand{\RR}{\mathbb{R}}

\newcommand{\FFF}{\mathcal{F}}

\begin{document}

\setlength{\parindent}{0cm}
\setlength{\parskip}{0.5cm}

\title[Hitting times for jump processes]{An elementary proof that the first hitting time of an $F_\sigma$ set by a
  jump process is a stopping time}

\author{Alexander Sokol}

\address{Alexander Sokol: Institute of Mathematics, University of
  Copenhagen, 2100 Copenhagen, Denmark}
\email{alexander@math.ku.dk}
\urladdr{http://www.math.ku.dk/$\sim$alexander}



\begin{abstract}
We give a short and elementary proof that the first hitting time of a
$F_\sigma$ set by the jump process of a c\`{a}dl\`{a}g adapted process is a
stopping time.
\end{abstract}

\maketitle

\noindent

\section{Introduction}

For a stochastic process $X$ and a subset $B$ of the real numbers, the
mapping $T$ defined by $T=\inf\{t\ge0|X_t \in B\}$ is called the first hitting time
of $B$ by $X$. In \cite{AS2013ELEM}, a short and elementary proof was given that the first
hitting time of an open set by the jump process of a c\`{a}dl\`{a}g
adapted process is a stopping time. A similar result is proved by
elementary means in \cite{MR1906715}, Proposition 1.3.14, where it is shown
that the first hitting time of $[c,\infty)$ for $c>0$ by the jump
process of a c\`{a}dl\`{a}g adapted process is a stopping time. Using
methods similar to both \cite{MR1906715} and \cite{AS2013ELEM}, we prove in this note that the hitting time of an
$F_\sigma$ set, meaning a countable union of closed sets, by the jump
process of a c\`{a}dl\`{a}g adapted process is a stopping time. As
open sets are $F_\sigma$ sets, this result covers both the case of
hitting an open and a closed set.

\section{Main result}

We assume given a filtered probability space $(\Omega,\mathcal{F},(\mathcal{F}_t),P)$
such that the filtration $(\mathcal{F}_t)_{t\ge0}$ is right-continuous in the
sense that $\mathcal{F}_t=\cap_{s>t}\mathcal{F}_s$ for all
$t\ge0$. Also, we use the convention that $X_{0-}=X_0$, so that there
is no jump at the timepoint zero.

\begin{theorem}
\label{theorem:FirstEntranceJumpOpen}
Let $X$ be a c\`{a}dl\`{a}g adapted process, and let $U$ be an $F_\sigma$ set in
$\mathbb{R}$. Define $T=\inf\{t\ge0|\varDelta X_t\in U\}$. Then $T$ is a
stopping time.
\end{theorem}
By the c\`{a}dl\`{a}g property of $X$, $\varDelta X$ is zero
everywhere except for on a countable set. Therefore, $T$ is identically
zero if $U$ contains zero, and so it is immediate that $T$ is a
stopping time in this case. We conclude that it suffices to prove the
result in the case where $U$ does not contain zero. Therefore, assume
that $U$ is an $F_\sigma$ set not containing zero. By right-continuity
of the filtration, it suffices to show $(T<t)\in\mathcal{F}_t$ for
$t>0$, see Theorem I.1 of \cite{MR2273672}. Fix $t>0$. Assume that $U=\cup_{n=1}^\infty F_n$, where $F_n$ is
closed. As $X_0-X_{0-}=0$ and $U$ does not contain
zero, we have
\begin{eqnarray}
  (T<t)&=&(\exists\; s\in(0,t):X_s-X_{s-}\in U)\notag\\
  &=&\cup_{s\in(0,t)}\cup_{n=1}^\infty (X_s-X_{s-}\in F_n)\notag\\
  &=&\cup_{n=1}^\infty \cup_{s\in(0,t)}(X_s-X_{s-}\in F_n)\notag\\
  &=&\cup_{n=1}^\infty(\exists\;s\in(0,t): X_s-X_{s-}\in F_n).
\end{eqnarray}
Thus, it suffices to show that $(\exists\;s\in(0,t): X_s-X_{s-}\in F)\in\FFF_t$ for all closed $F$. Assume given
such a closed set $F$. We claim that
\begin{eqnarray}
  &&(\exists\; s\in(0,t):X_s-X_{s-}\in F)
 =\cap_{n=1}^\infty \cup_{(p,q)\in\Theta_n} (X_q-X_p\in F_n)\;\label{eq:Equality},
\end{eqnarray}
where $F_n = \{x\in \mathbb{R}\mid \;\exists\;y\in F:|x-y|\le 1/n\}$ and $\Theta_n$ is the subset of $\mathbb{Q}^2$ defined by $\Theta_n =
\{(p,q)\in\mathbb{Q}^2|0<p<q<t,|p-q|\le1/n\}$.

To prove this, we first consider the inclusion towards the right. Assume
that there is $0<s<t$ such that $X_s-X_{s-}\in F$. Fix $n\ge1$. By
the path properties of $X$, we obtain that for $p,q\in\RR$ with $0<p<s<q<t$ and $p$ and
$q$ close enough to $s$, $|X_q-X_s|\le 1/2n$ and
$|X_p-X_{s-}|\le 1/2n$, yielding $|(X_q-X_p)-(X_s-X_{s-})|\le
1/n$ and thus $X_q-X_p\in F_n$. By picking $p$ and $q$ in $\QQ$ close
enough to $s$, we obtain $(p,q)\in\Theta_n$ as well. This proves the
inclusion towards the right.

Next, consider the inclusion towards the
left. Assume that for all $n\ge1$, there is $(p_n,q_n)\in\Theta_n$ such
that $X_{q_n}-X_{p_n}\in F_n$. We then also have $\lim_n
|p_n-q_n|=0$. By taking two consecutive subsequences
and relabeling, we may assume that in addition to having $\lim_n
|p_n-q_n|=0$ and $0<p_n<q_n<t$, both $p_n$ and $q_n$ are
monotone. As $(F_n)$ is decreasing, we then also obtain
$X_{q_n}-X_{p_n}\in F_n$ for all $n\ge1$. As $p_n$ and $q_n$ are
bounded and monotone, they are convergent, and as $\lim_n|q_n-p_n|=0$,
it follows that the limit $s$ is the same for both $q_n$ and $p_n$.

We wish to argue that $0<s<t$, that $X_{s-}=\lim_nX_{p_n}$ and that $X_s=\lim_n
X_{q_n}$. First note that as both $(p_n)$ and $(q_n)$ are monotone, the limits
$\lim_nX_{p_n}$ and $\lim_nX_{q_n}$ exist and are either equal to $X_s$ or
$X_{s-}$. As $X_{q_n}-X_{p_n}\in F_n$, we obtain
\begin{displaymath}
  \lim_nX_{q_n}-\lim_nX_{p_n}
  =\lim_n X_{q_n}-X_{p_n}\in \cap_{n=1}^\infty F_n=F,
\end{displaymath}
where the final equality follows as $F$ is closed. As $F$ does not contain zero, we conclude
$\lim_nX_{q_n}-\lim_nX_{p_n}\neq0$. From this, we immediately obtain $0<s<t$, as if $s=0$, we
would obtain that both $\lim_n X_{q_n}$ and $\lim_nX_{p_n}$ were equal
to $X_s$, and if $s=t$, both $\lim_n X_{q_n}$ and $\lim_nX_{p_n}$
would be equal to $X_{s-}$, in both cases yielding a contradiction. Also, we cannot have that both limits are
$X_s$ or that both limits are $X_{s-}$, and so only two cases are
possible, namely that $X_s=\lim_nX_{q_n}$ and $X_{s-}=\lim_nX_{p_n}$ or that
$X_s=\lim_nX_{p_n}$ and $X_{s-}=\lim_nX_{q_n}$. We wish to argue that
the former holds. If $X_s=X_{s-}$, this is trivially the
case. Assume that $X_s\neq X_{s-}$ and that $X_s=\lim_nX_{p_n}$ and
$X_{s-}=\lim_nX_{q_n}$. If $q_n\ge s$ eventually or $p_n<s$ eventually,
we obtain $X_s=X_{s-}$, a contradiction. Therefore, $q_n<s$ infinitely often and $p_n\ge s$ infinitely
often. By monotonicity, $q_n<s$ and $p_n\ge s$ eventually, a
contradiction with $p_n<q_n$. We conclude $X_s=\lim_nX_{q_n}$ and
$X_{s-}=\lim_nX_{p_n}$, as desired.

From this, we conclude $X_s-X_{s-}=\lim_n X_{q_n}-\lim_n X_{p_n}\in
F$. This proves the existence of $s\in(0,t)$ such that $X_s-X_{s-}\in
F$, and so proves the inclusion towards the right.

We have now shown (\ref{eq:Equality}). Now, as $X_s$ is
$\mathcal{F}_t$ measurable for all $0\le s\le t$, the set
$\cap_{n=1}^\infty \cup_{(p,q)\in\Theta_n} (X_q-X_p\in F_n)$ is
$\mathcal{F}_t$ measurable as well. We conclude that $(T<t)\in\mathcal{F}_t$ and so $T$ is a stopping time.

\bibliographystyle{amsplain}

\bibliography{full}

\end{document}